# ON THE TOPOLOGY OF COMPACT STEIN SURFACES

SELMAN AKBULUT AND BURAK OZBAGCI

ABSTRACT. In this paper we obtain the following results: (1) Any compact Stein surface with boundary embeds naturally into a symplectic Lefschetz fibration over $S^2$. (2) There exists a minimal elliptic fibration over $D^2$, which is not Stein. (3) The circle bundle over a genus $n \geq 2$ surface with euler number $e = -1$ admits at least $n+1$ mutually non-homeomorphic simply-connected Stein fillings. (4) Any surface bundle over $S^1$, whose fiber is a closed surface of genus $n \geq 1$ can be embedded into a closed symplectic 4-manifold, splitting the symplectic 4-manifold into two pieces both of which have positive $b_2^+$. (5) Every closed, oriented connected 3-manifold has a weakly symplectically fillable double cover, branched along a 2-component link.

## 0. INTRODUCTION

In [AO] (see also [LP]), the authors proved that every compact Stein surface admits a PALF (positive allowable Lefschetz fibration over $D^2$ with bounded regular fibers) and conversely every PALF is Stein. In this paper we first prove that any compact Stein surface with boundary embeds naturally into a symplectic Lefschetz fibration over $S^2$. In particular, this shows that we can embed a Stein surface into a closed symplectic 4-manifold.

Next we show that the result in [AO] does not necessarily hold if the fiber of the Lefschetz fibration is closed, by constructing non-Stein minimal elliptic fibrations (with closed regular fibers) over $D^2$. Minimality of our examples is important, otherwise one can easily find non-Stein examples by blowing-up an elliptic (or in general a Lefschetz) fibration over $D^2$.

We also prove that for every integer $n \geq 2$, there exits an irreducible 3-manifold $M_n$ with at least $n+1$ mutually non-homeomorphic simply-connected Stein fillings. We will identify $M_n$ as the circle bundle over a genus $n$ surface with euler number $e = -1$.

Moreover we show that any surface bundle over $S^1$, whose fiber is a closed surface of genus $n \geq 1$ can be embedded into a closed symplectic 4-manifold, splitting the symplectic 4-manifold into two pieces both of which have positive $b_2^+$.

Finally we prove that every closed, oriented connected 3-manifold has a weakly symplectically fillable double cover, branched along a 2-component link. This result is interesting in our point of view, since not every 3-manifold is weakly symplectically fillable as shown in [L]. We advise the reader to turn to [EH] for notions like weak/strong symplectic filling, Stein filling, etc.

## 1. Natural symplectic compactifications of Stein surfaces

**Theorem 1.** *Any compact Stein surface with boundary embeds naturally into a symplectic Lefschetz fibration over $S^2$. In particular, any compact Stein surface with boundary embeds into a closed symplectic 4-manifold. Conversely given a symplectic Lefschetz fibration over $S^2$ with a section, then if we remove a neighborhood of this section union a regular fiber we get a Stein surface.*

*Proof.* We know that a compact Stein surface $X$ with boundary admits a PALF. We may also assume that the regular fiber $F$ has only one boundary component. The fibration induces an open book decomposition of $\partial X$ with connected binding $\partial F$. First we enlarge $X$ to $X'$ by attaching a 2-handle along $\partial F$ with 0-framing. Note that $\partial X'$ is an $\hat{F}$-bundle over $S^1$, where $\hat{F}$ denotes the closed surface obtained by capping off the surface $F$ by gluing a 2-disk along its boundary. Also $X'$ admits a positive Lefschetz fibration over $D^2$ with regular fiber $\hat{F}$. Let $\mathrm{Map}(\hat{F})$ denote the mapping class group of the closed surface $\hat{F}$. The second author learned the proof of the next lemma from Ivan Smith.

**Lemma 2.** *Any element in $\mathrm{Map}(\hat{F})$ can be expressed as a product of nonseparating positive Dehn twists.*

*Proof.* Let the curves $A_i, B_i$ on a genus $n \geq 1$ surface $\hat{F}$ be drawn as in Figure 1 and write $a_i$ and $b_i$ for the positive Dehn twists about $A_i$ and $B_i$, respectively.
The following is a standard word in the mapping class group $\mathrm{Map}(\hat{F})$. (cf. [B]).

$$1 = (a_1 b_1 a_2 b_2 ... a_n b_n)^{4n+2}.$$

Hence we conclude that $a_1^{-1}$ is a product of nonseparating positive Dehn twists. Therefore any negative nonseparating Dehn twist is a product of nonseparating positive Dehn twists since any two negative nonseparating Dehn twists are equivalent. This finishes the proof of the lemma combined with the fact that $\mathrm{Map}(\hat{F})$ is generated by (positive and negative) nonseparating Dehn twists.
□

We use Lemma 2 to extend $X'$ into a positive Lefschetz fibration $\tilde{X}$ over $S^2$ with regular fiber $\hat{F}$ as follows. Let $c_1 c_2 ... c_k$ be the global monodromy of the PALF on $X$, where $c_i$ denotes the positive Dehn twist along a simple closed curve $C_i$ on $F$. Then this product (after capping off the boundary component) can be viewed as a product in $\text{Map}(\hat{F})$. We clearly have

$$c_1 c_2 ... c_k c_k^{-1} c_{k-1}^{-1} ... c_1^{-1} = 1.$$

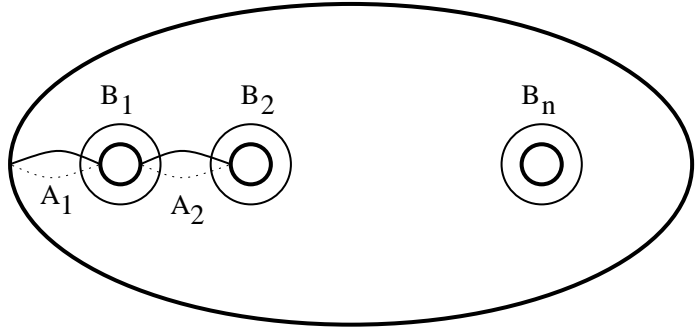

Figure 1. A surface of genus $n \geq 1$

By Lemma 2 we can replace every negative twist by a product of positive twists to obtain $\tilde{X}$. $\tilde{X}$ admits a symplectic structure with symplectic regular fibers by a Theorem of Gompf [GS]. Consequently, $X$ is emdedded naturally into a closed symplectic 4-manifold $\tilde{X}$.

Conversely suppose that $\tilde{X}$ is a symplectic Lefschetz fibration over $S^2$ with a section. Then if we remove a neighborhood of this section union a regular fiber we get a Stein surface. This is clear from our description of Stein surfaces as PALF's.

□

*Remark* 1. A version of Theorem 1 was proved in [LM]. Our proof, however, is more topological and should be viewed as the symplectic version of the standard procedure of embedding a compact manifold into a closed manifold by doubling process. In fact, this theorem gives a procedure of associating to any PALF a "maximal" positive Lefschetz fibration over $D^2$ with closed fibers, which is then capped off by $\hat{F} \times D^2$ to get a symplectic Lefschetz fibration over $S^2$. In particular, we have a map

$$\{\text{PALF's}\} \to \{\text{Symplectic 4-manifolds}\}.$$

Moreover, since a compact Stein surface $X$ admits infinitely many nonequivalent PALF's (cf. [AO]), we can embed $X$ into a symplectic Lefschetz fibration of arbitrarily large genus over $S^2$. Any fibered knot $K$ in $S^3$ with positive allowable monodromy induces a $(PALF)_K$ on $B^4$ and hence the association above gives a closed symplectic 4-manifold. For example,

$$(PALF)_{trefoil} \to K3 \text{ surface}$$

To see this implication consider the the PALF induced by the trefoil knot as indicated in Figure 2. Its monodromy $ab$ is given by Dehn twists along the two generators of the fiber which is a punctured torus.

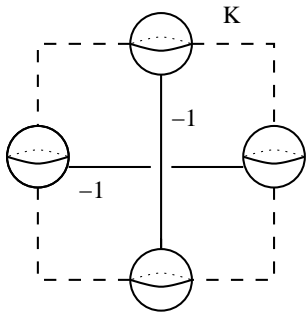

FIGURE 2. PALF induced by trefoil

Since $(ab)^6 = 1$, we have $a^{-1} = b(ab)^5$ and $b^{-1} = (ab)^5 a$, hence by our algorithm $abb^{-1}a^{-1} = (ab)^{12}$. To get the symplectic completion we need to attach 0 framed 2−handle to the binding (trefoil knot) and 22 more 2−handles with −1 framing each, as indicated in Figure 3. We then cap it off with $D^2 \times T^2$ at the end, to get a closed manifold, which corresponds to attaching the indicated −2 framed handle in Figure 3 and two 3−handles (which we don't need to indicate). This gives K3 surface (compare [HKK]) which is usually denoted by E(2).

Notice, that since $(ab)^{-1} = (ab)^5$, rather then using the algorithm we can write in a shorter way $abb^{-1}a^{-1} = (ab)^6$, and end up constructing a smaller symplectic completion by attaching ten 2-handles to Figure 2 instead, as indicated in Figure 4. This time the last 2-handle corresponding to $D^2 \times T^2$ has to be attached by −1 framing (see discussion in [HKK] about determining framing of this handle). This gives us the half Kummer surface $E(1)$.

*Remark* 2. The last statement in Theorem 1 is implicit in Donaldson's work ([D]). Our proof, however, is purely topological. Combining with the existence of Lefschetz pencils (cf. [D]) on symplectic 4-manifolds one concludes that every closed symplectic 4-manifold, possibly after blowing up at some points, contains a "big" Stein piece.

*Remark* 3. In Lemma 2 we proved that $\text{Map}(\hat{F})$ is generated by nonseparating positive Dehn twists. This is not true for the mapping class group of a surface with boundary. Otherwise any 3-manifold would be Stein fillable (cf. [AO]), which is shown to be false, for example, in Theorem 7.

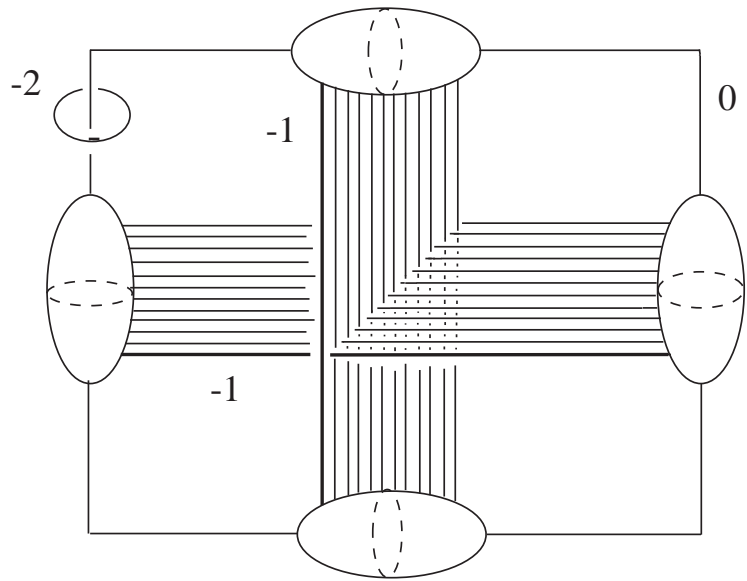

Figure 3. E(2)

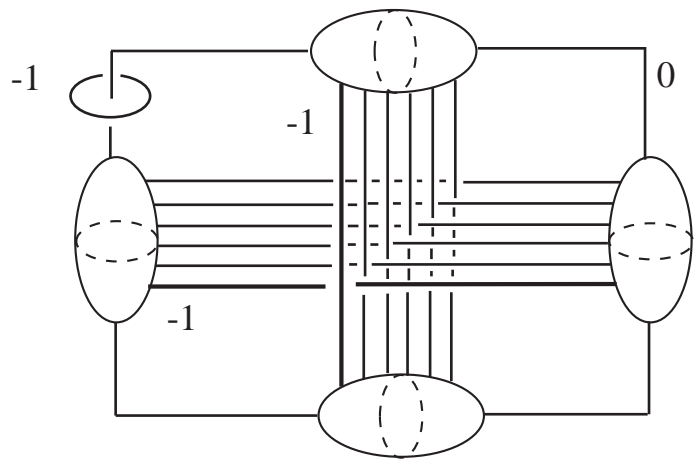

Figure 4. E(1)

## 2. Non-Stein minimal elliptic fibrations

Let $Y(e,n)$ denote the circle bundle over a genus $n \geq 1$ surface with euler number $e$. First we observe that $Y(1,1)$ is also a torus bundle over $S^1$ with monodromy a single positive Dehn twist along a simple closed curve on the torus $T^2$. (A proof of this is given in the appendix). Next we show that $Y(1,1)$ is Stein fillable. In fact a regular neighborhood of so called fishtail fiber in an elliptic fibration is a Stein filling of $Y(1,1)$ as shown in Figure 5.

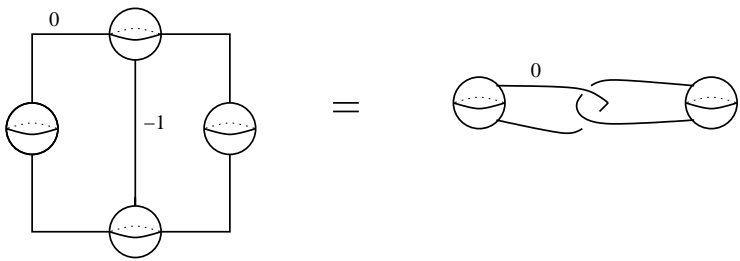

Figure 5. Regular neigborhood of a fishtail fiber

Now consider the elliptic surface $E(3) = E(1) \natural E(1) \natural E(1)$ which is a 3-fold fiber sum of the rational elliptic surface $E(1) = \mathbb{C}P^2 \# 9\overline{\mathbb{C}P^2}$. (Here $\natural$ denotes the fiber sum and $\#$ denotes the connected sum.) The global monodromy of the elliptic fibration of $E(3)$ over $S^2$ is given by $(ab)^{18}$, where $a$ and $b$ denote positive Dehn twists along simple closed curves $A$ and $B$, respectively as shown in Figure 6. Let $X$ denote the elliptic fibration over $D^2$ with monodromy $(ab)^6 a$. Then $E(3)$ splits as a union of $X = E(1) \natural Z_1$ and $Z_2 \natural E(1)$ such that $\partial X = Y(1,1)$.

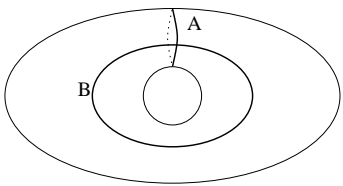

Figure 6. Torus

**Theorem 3.** *X is a minimal, symplectic and simply-connected elliptic fibration over $D^2$, which does not admit any Stein structure.*

*Proof.* $X$ is minimal since $E(3)$ is minimal. $X$ is simply-connected since $X = E(1) \natural Z_1$ and $E(1)$ is simply-connected. $X$ admits a symplectic structure by a theorem of Gompf [GS].

Suppose that $X$ is Stein. We can embed $X$ into a complex surface $S$ of general type by a theorem of Lisca and Matic [LM]. We will show that we can assume $b_2^+(S-X) > 0$. Suppose that $b_2^+(S-X) = 0$. First we extend $X$ into $X'$ by attaching a 2-handle with framing $tb(K) - 1$ along a Legendrian knot $K$ satisfying $tb(K) > 1$ and contained in a standard 3-ball $D^3$ in $\partial X$. (Here $tb$ denotes the Thurston-Bennequin invariant). Because of the framing of the 2-handle, $X'$ is also a Stein surface. So we can embed $X'$ into a complex surface $S'$ of general type. (cf. Figure 7). Thus $b_2^+(S' - X) > 0$, since we created a second homology class with positive self-intersection in $S' - X$. Also note that $b_2^+(X) > 0$.

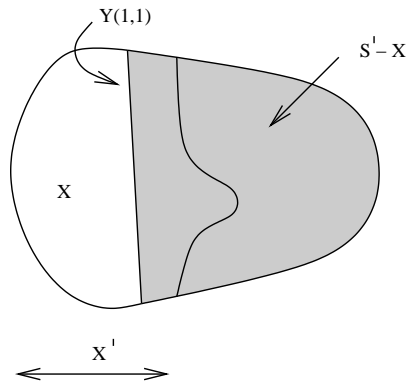

FIGURE 7.

This gives a contradiction to the following result of Ozsvath and Szabo, hence showing that $X$ is not Stein.

**Theorem 4.** [OzS] *If $S$ is a complex surface of general type, and $|e| \geq 2n-1$, then $S$ admits no splitting along an embedded copy of $Y = Y(e,n)$ of the form $S = X_1 \#_Y X_2$ with $b_2^+(X_1), b_2^+(X_2) > 0$.*

□

*Remark* 4. We can generalize Theorem 3 as follows. Define $X_k = E(k) \natural Z_1$ for $k \geq 1$. Then for each $k \geq 1$, $X_k$ is a minimal, symplectic and simply-connected elliptic fibration over $D^2$, which does not admit any Stein structure. We can also generalize Theorem 3 in a different direction. One can easily find non-Stein *achiral* minimal elliptic fibrations over $D^2$.

**Theorem 5.** *Suppose that $\tilde{X}$ is a Stein filling of some $Y(e,n)$ with $|e| \geq 2n - 1$, then $b_2^+(\tilde{X}) = 0$.*

*Proof.* We use the same trick as in the proof of Theorem 3. □

The following result was proved in [GS]. We give a different proof as an application of Theorem 5.

**Corollary 6.** *Let $X(e,n)$ denote the disk bundle over a genus $n$ surface with euler number $e$. If $e \geq 2n - 1$ then $X(e,n)$ is not Stein.*

*Proof.* First note that $\partial X(e,n) = Y(e,n)$. Suppose that $e \geq 2n - 1$. Then we can apply Theorem 5 to conclude that $X(e,n)$ is not Stein since $b_2^+(X(e,n)) = 1$. □

*Remark* 5. [GS] The result in Corollary 6 is optimal, i.e., $X(e,n)$ is realized as a Stein surface iff $e < 2n - 1$. $Y(e,n) = \partial X(e,n)$, however, is Stein fillable for all $e, n$.

It is known that the Poincare homology sphere with reversed orientation is not weakly symplectically semi-fillable [L]. We prove the following weaker result by our technique.

**Theorem 7.** *Poincare homology sphere with reversed orientation is not Stein fillable.*

*Proof.* Let Y denote the Poincare homology sphere oriented as the boundary of the $-E_8$ plumbing as shown in Figure 8.

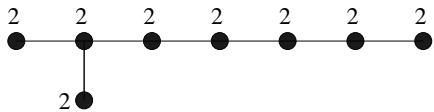

FIGURE 8. $-E_8$ plumbing

Note that $Y$ has a metric with positive scalar curvature. Suppose $X$ is a Stein filling of $Y$. Embed $X$ into a complex surface $S$ of general type with $b_2^+(S - X) > 0$ as in the proof of Theorem 3. Since $S$ has non-vanishing Seiberg-Witten invariants and $Y$ has a metric with positive scalar curvature, we have $b_2^+(X) = 0$. (see [OhO] for a proof). Hence $X \cup (E_8)$ is a closed, smooth and negative definite 4-manifold, which can not exist by a theorem of Donaldson. □

## 3. Stein fillings of some circle bundles

In this section we show that for every integer $n \geq 2$, there exits an irreducible 3-manifold $M_n$ with at least $n+1$ mutually non-homeomorphic simply-connected Stein fillings. We will also identify $M_n$ as the circle bundle $Y(-1, n)$ over a genus $n$ surface with euler number $e = -1$.

Let $t_\alpha$ denote a positive Dehn twist about a simple closed curve $\alpha$ on an oriented surface $F$. The following result is standard. (cf. [B]).

**Lemma 8.** *For any two simple closed curves $\alpha$ and $\beta$ on $F$ we have $t_\beta t_\alpha = t_\alpha t_{t_\beta(\alpha)}$.*

Let the curves $A_i, B_i, D_2, E_2$ on a genus $n \geq 2$ surface $F$ with one boundary component be drawn as in Figure 9 and write lower-case letters $a_i$ etc. for the positive Dehn twist about $A_i$ etc. Let $\delta$ denote a positive Dehn twist about a curve parallel to the boundary.

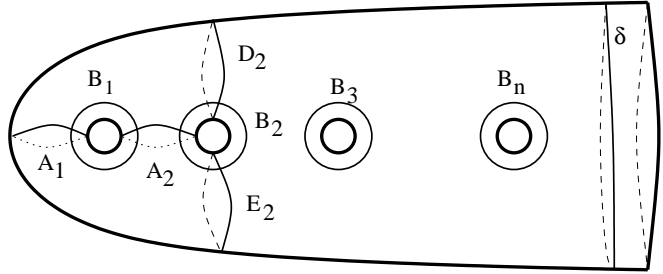

FIGURE 9. Genus $n \geq 2$ surface $F$ with boundary

**Proposition 9.**
$$(a_1 b_1 a_2 b_2 ... a_n b_n)^4 = (a_1 b_1 a_2)^4 \psi,$$
where $\psi$ is a product of $8n - 12$ nonseparating positive Dehn twists.

*Proof.* We apply Lemma 8 as many times as needed to obtain $(a_1 b_1 a_2)^4$ at the beginning of the given product. $\square$

**Lemma 10.** [LP] *We have $\delta = (a_1 b_1 a_2 b_2 ... a_n b_n)^{4n+2}$ in $Map(F, \partial F)$.*

**Theorem 11.** *For every integer $n \geq 2$, there exits an irreducible 3-manifold $M_n$ with at least $n+1$ mutually non-homeomorphic simply-connected Stein fillings.*

*Proof.* Fix an integer $n \geq 2$. Define
$$F_\delta = (F \times [0,1])/(\delta(x), 0) \sim (x, 1)$$
where $\delta \in Map(F, \partial F)$ is the map in Lemma 10. Consider the closed oriented 3-manifold
$$M_n = (S^1 \times D^2) \cup F_\delta.$$
In other words, $M_n$ has a positive open book decomposition with binding $\partial F$, page $F$ and monodromy $\delta$ which is expressed as a product of nonseparating positive Dehn twists. By the results in [AO], $M_n$ admits a Stein filling $X_0$ which is a PALF (positive allowable Lefschetz fibration over $D^2$ with bounded fibers) of genus $n$ with $2n(4n+2)$ singular fibers. We can calculate the Euler characteristic $\chi(X_0)$ as follows:

$$\chi(X_0) = 2 - 2n - 1 + 2n(4n+2) = 8n^2 + 2n + 1.$$

The following is a standard relation in the mapping class group. (cf. [B]).

$$(a_1 b_1 a_2)^4 = d_2 e_2.$$

Using this relation and Proposition 9, we can write the map $\delta$ as follows:

$$\delta = (a_1 b_1 a_2 b_2 ... a_n b_n)^{4n+2}$$
$$= (a_1 b_1 a_2 b_2 ... a_n b_n)^4 (a_1 b_1 a_2 b_2 ... a_n b_n)^{4(n-1)+2}$$
$$= (a_1 b_1 a_2)^4 \psi (a_1 b_1 a_2 b_2 ... a_n b_n)^{4(n-1)+2}$$
$$= d_2 e_2 \psi (a_1 b_1 a_2 b_2 ... a_n b_n)^{4(n-1)+2}$$

This last product gives another Stein filling $X_1$ of $M_n$ as a PALF with 10 less singular fibers. Thus

$$\chi(X_1) = \chi(X_0) - 10.$$

We can iterate the substitution above $n$-times to get $n+1$ mutually non-homeomorphic Stein fillings $X_0, X_1, ..., X_n$ of $M_n$ such that

$$\chi(X_i) = \chi(X_{i-1}) - 10.$$

It is easy to see that all the Stein fillings are simply-connected. Next we show that $M_n = Y(-1, n)$, which implies in particular, that $M_n$ is irreducible. □

**Lemma 12.** [St] *We can identify $M_n$ as the circle bundle $Y(-1, n)$ over a genus $n$ surface with euler number $e = -1$.*

*Proof.* Fix an integer $n \geq 2$. Let $F$ denote a genus $n$ surface with one boundary component and let $\hat{F}$ denote the closed surface obtained by capping off the surface $F$ by gluing a 2-disk along its boundary. Then there is a natural map

$$Map(F, \partial F) \to Map(\hat{F}).$$

So the relation in Lemma 10 induces a word in $Map(\hat{F})$. Consequently this gives a positive Lefschetz fibration over $S^2$ with a section a sphere of square $-1$ and regular fiber $\hat{F}$. Now consider a neighborhood $U_n$ of a regular fiber union this section. First we observe that $\partial U_n = \overline{M_n}$. Moreover $U_n$ is obtained by plumbing a disk bundle over $\hat{F}$ and a disk bundle over $S^2$. We can blow down the $-1$ sphere to get a disk bundle over $\hat{F}$ with euler number $+1$. We prove our result by reversing the orientations. □

*Remark* 6. In [G], Gompf proves that $Y(-1, n)$ has at least $n$ Stein fillable contact structures. The Stein fillings are mutually diffeomorphic but they admit different complex structures distinguished by their first Chern classes.

4. SPLITTINGS OF SYMPLECTIC 4-MANIFOLDS ALONG SURFACE BUNDLES OVER $S^1$

**Proposition 13.** *Any closed surface bundle over $S^1$ can be embedded into a closed symplectic 4-manifold, splitting the symplectic 4-manifold into two pieces both of which have positive $b_2^+$.*

*Proof.* Let $\Sigma$ denote a closed, connected and oriented surface of genus $n$. Let $B_\phi$ denote the $\Sigma$-bundle over $S^1$ with monodromy $\phi \in \mathrm{Map}(\Sigma)$. We can express $\phi$ as a product $c_1 c_2 ... c_k$ of (nonseparating) positive Dehn twists by Lemma 2. Using this product we can fill in the $\Sigma$-bundle over $S^1$ with a positive Lefschetz fibration $X_1$ over $D^2$ with regular fiber $\Sigma$. We trivially have

$$c_1 c_2 ... c_k c_k^{-1} c_{k-1}^{-1} ... c_1^{-1} = 1.$$

Now replace each negative twist in this word by a product of (nonseparating) positive Dehn twists again by Lemma 2. Hence we get a positive Lefschetz fibration $X$ over $S^2$ which is a union of $X_1$ and $X_2 = X - X_1$. Suppose that $b_2^+(X_i) = 0$ for some $i = 1, 2$. (Otherwise the theorem is proved). Let $G(n)$ denote the genus $n$ Lefschetz fibration over $S^2$ given by the word

$$1 = (a_1 b_1 a_2 b_2 ... a_n b_n)^{4n+2}.$$

These higher genus Lefschetz fibrations can be considered as the generalization of the elliptic fibration on $E(1)$. We have $b_2^+(G(n)) > 0$ since $G(n)$ has a symplectic structure by a theorem of Gompf [GS]. Define

$$X_i' = X_i \sharp G(n)$$

where $\sharp$ denotes the fiber sum along a regular fiber. It is easy to see that $b_2^+(X_i') > 0$ for $i = 1, 2$. Hence $G(n) \sharp X \sharp G(n)$ is a closed symplectic 4-manifold which is a union of $X_1'$ and $X_2'$ (with $b_2^+(X_i') > 0$), glued along the $\Sigma$-bundle $B_\phi$. □

*Remark* 7. We were pointed out by I. Smith that the first part of our proof (not including the $b_2^+$ arguement) had appeared in [Sm]. He actually shows that the Lefschetz fibration $X_1$ is a weak symplectic filling for the contact structure arising from a $C^0$ perturbation of the obvious foliation on $\partial X_1$.

In Theorem 1, we described a procedure to obtain a symplectically fillable 3-manifold by surgery from an open book decomposition. In the following, we use branched coverings to obtain a similar result.

**Theorem 14.** *Every closed, oriented connected 3-manifold has a weakly symplectically fillable double cover, branched along a 2-component link.*

*Proof.* We use the following result of Montesinos [M] (see also [S]).

**Lemma 15.** [M] *Every closed, oriented connected 3-manifold $M$ contains a two component link $L$, such that there is a 2-fold covering of $M$ branched over $L$ which is a surface bundle over $S^1$.*

*Proof.* It is well known that every closed, oriented connected 3-manifold $M$ has an open book decomposition with connected binding. Let $F$ denote the page and $\phi$ denote the monodromy of an open book on $M$. We use $2F$ to denote the double of $F$, where two copies of $F$ are glued along their boundary. Define
$$N = 2F \times [-1, 1]/ \ (\phi\#\phi^{-1}(x), -1) \sim (x, 1)$$
So, $N$ is a $2F$-bundle over $S^1$ with monodromy $\phi\#\phi^{-1}$, where $\phi$ acts on one copy of $F$ and $\phi^{-1}$ acts on the other copy. Let $u$ be the involution on $N$ given by
$$u(x, t) = (ix, -t),$$
where $i$ denotes the involution on $2F$ interchanging the two copies of $F$. Then it is easy to see that $N/u = M$. Moreover the fixed point set of the action of $u$ on $N$ is a disjoint union of two circles. □

Since every 3-manifold $M$ is double branched covered by a surface bundle $N$, and every surface bundle over $S^1$ is weakly symplectically fillable as in the proof of Proposition 13, we conclude that every closed oriented connected 3-manifold has a weakly symplectically fillable double cover, branched along a 2-component link. □

## 5. Appendix

We prove that the circle bundle over torus with euler number $e = 1$ is diffeomorphic to a torus bundle over circle, whose monodromy is a single positive Dehn twist.

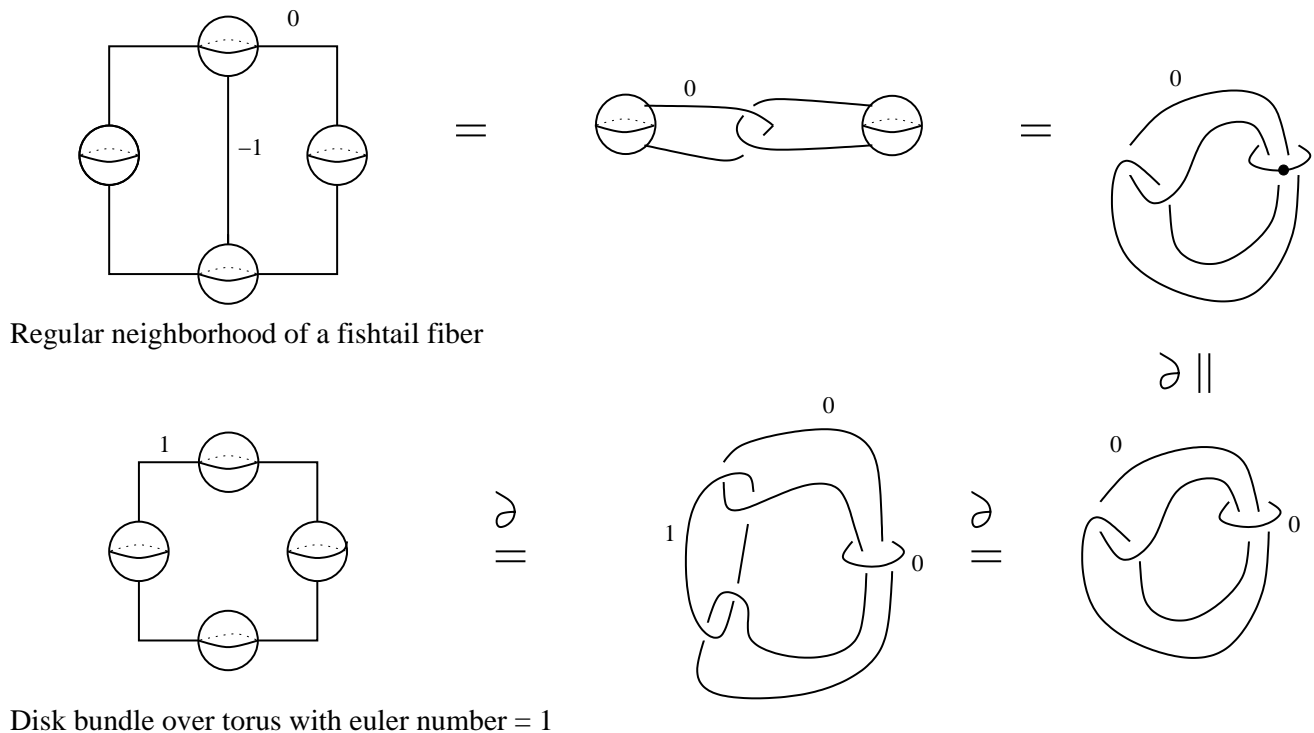

Regular neighborhood of a fishtail fiber

Disk bundle over torus with euler number = 1

Figure 10.

Department of Mathematics, Michigan State University, MI, 48824
*E-mail address*: akbulut@math.msu.edu and bozbagci@math.msu.edu